\documentclass[11pt]{amsart}

\usepackage{amssymb,xfrac}
\usepackage{enumerate}
\usepackage{tikz}
\usetikzlibrary{arrows}
\usepackage{graphicx}
\usepackage{hyperref}
\usepackage{ulem}
\usepackage{tikzit}

\tikzstyle{red dot}=[fill=red, draw=black, shape=circle]
\tikzstyle{green dot}=[fill={rgb,255: red,0; green,174; blue,0}, draw=black, shape=circle]
\tikzstyle{black dot}=[fill=black, draw=black, shape=circle]
\tikzstyle{blackoutline}=[fill=white, draw=black, shape=circle, minimum size=2mm]

\tikzstyle{bold edge}=[-, fill=black, tikzit fill=black, thick]
\tikzstyle{grey edge 1}=[-, fill={rgb,255: red,128; green,128; blue,128}, tikzit fill={rgb,255: red,128; green,128; blue,128}, draw={rgb,255: red,128; green,128; blue,128}]
\tikzstyle{black dash 1}=[-, thick, dashed, dash pattern=on 4mm off 2mm]
\tikzstyle{little dots}=[-, thick, dashed, dash pattern=on 1mm off 1mm]

\newtheorem{thm}{Theorem}[section]

\newtheorem{lem}[thm]{Lemma}
\newtheorem{cor}[thm]{Corollary}

\newtheorem{conj}[thm]{Conjecture}

\def\ex{{\operatorname{ex}}}
\def\exr{{\operatorname{ex}^\star}}

\def\spec{{\operatorname{Spec}}}
\def\ds{{DS_{r,s}}}

\theoremstyle{definition}

\theoremstyle{remark}

\newtheorem{remark*}{Remark}

\numberwithin{equation}{section}

\title{Rainbow Tur\'{a}n Methods for Trees}
\author{Vic Bednar
\and Neal Bushaw}
\thanks{
Department of Mathematics and Applied Mathematics,
Virginia Commonwealth University, USA
}

\begin{document}

\maketitle

\begin{abstract}
The rainbow Tur\'{a}n number, a natural extension of the well studied traditional Tur\'{a}n number, was introduced in 2007 by Keevash, Mubayi, Sudakov and Verstra\"te.  The rainbow Tur\'{a}n number of a graph $H$, $\exr(n,H)$, is the largest number of edges for an $n$ vertex graph $G$ which can be properly edge colored with no rainbow $H$ subgraph. 

We explore the reduction method for finding upper bounds on rainbow Tur\'{a}n numbers, and use this to inform results for the rainbow Tur\'{a}n numbers of double stars, caterpillars, and perfect binary trees. In addition, we define $k$-unique colorings and the related $k$-unique Tur\'{a}n numbers. We provide preliminary results on this new variant on the classic problem.
\end{abstract}

\section{Introduction}
The forbidden subgraph problem is perhaps the prototypical question in extremal graph theory.  How do we maximize the number of edges in an $n$-vertex graph, while forbidding some fixed subgraph? 

The earliest famous instance of this sort of problem dates to 1907; now known as Mantel's Theorem, this gives the largest possible number of edges in an $n$ vertex graph that does not contain any triangles \cite{mantel}\footnote{It is worth noting that although Mantel is typically given credit for this reference, the history is somewhat more complicated. Mantel submitted the problem to {\emph{Nieuwe Opgaven}}, where it was distributed loose-leaf to the members of the Royal Dutch Mathematical Society.  Eventually, this was reprinted in a digest version, along with a solution due to Wythoff (although it is further noted in this version that solutions were also submitted by H. Gouwentak, W. Mantel, J. Teixeira de Mattos, and Dr. F. Schuh, but only Wythoff's proof appears).}. In the 1940's Tur\'{a}n generalized Mantel's result to forbidding any particular complete graph, rather than just the triangle. In order to state such results compactly, we define the Tur\'{a}n number of a graph, $\ex(n, F)$\footnote{From this point forward we will refer to this as the traditional Tur\'{a}n number.}, to be the largest number of edges on $n$ vertices which is $F$-free for some subgraph $F$. In this context, we refer to $F$ as the forbidden subgraph.

The Erd\H{o}s-Stone theorem states that $\ex(n,F) = \frac{1}{\chi(F)-1}\binom{n}{2} + o(n^2)$, where $\chi(F)$ is the chromatic number of the forbidden subgraph \cite{E-Stone}. As the theorem gives degenerate results for graphs whose chromatic number is two, much of the research into traditional Tur\'{a}n numbers is centered around bipartite graphs. While many results address specific bipartite graphs or specific families of bipartite graphs, there is still little known for bipartite graphs as a whole (for an exhaustive survey, see \cite{FS}).

In this paper we consider two closely-related variants of the traditional Tur\'{a}n problem. The first, whose systemic study was formalized in 2007 by Keevash, Mubayi, Sudakov, and Verstra\"{e}te, adds in a proper edge coloring of the host graph and requires that the forbidden subgraph be ``rainbow'' -- every edge colored a distinct color. The {\bf{rainbow Tur\'{a}n number}} of a graph $F$, denoted $\exr(n,F)$, is the largest possible number of edges among those $n$ vertex graphs which can be properly edge-colored in a way that contains no rainbow $F$ subgraph. There are two important subtleties in this seemingly simple definition. First, there must exist a graph with $n$ vertices and $\exr(n, F)$ edges which admits at least one proper edge coloring that does not contain a rainbow $F$ subgraph. Additionally, every proper edge coloring of every graph with $n$ vertices and at least $\exr(n,F) +1$ edges must contain a rainbow $F$ subgraph.

The second variation (introduced in this paper) considers the case where a specific number of edges in $F$ must be assigned distinct colors. Given graphs $G, F$, a natural number $k$, and a proper edge coloring $\phi~:~E(G) \rightarrow \mathbb{N}$, we say that $G$ contains an {\bf{exactly $k$-unique copy of $F$}} if there are exactly $k$ edges whose color appears exactly once on $F$. Similarly, $G$ contains a {\bf{$k$-unique copy of $F$}} when there are \textit{at least} $k$ edges whose color appears exactly once on $F$. We define a new paramaterized family of extremal functions, which run from the traditional Tur\'{a}n number to the rainbow Tur\'{a}n number. Given a forbidden graph $F$ and a natural number $k$, the {\bf{$k$-unique Tur\'{a}n number}} $\ex_{k}(n, F)$ is the largest number of edges in an $n$ vertex graph which has some proper edge coloring that contains no $k$-unique $F$.\footnote{We note that throughout, we will be interested in non-exact colorings -- that is, $k$-unique as opposed to exactly $k$-unique. A larger number of uniquely colored edges will be good for us.  We also emphasize that throughout, we mean \emph{edge}-colorings when we refer to colorings.}


In this paper we consider the rainbow Tur\'{a}n problem and the $k$-unique Tur\'{a}n problem for families of small trees. In particular, we consider the family of double stars. For any $r,s \in \mathbb{N}$, the double star $\ds$ is the tree formed by taking a single edge $yx$ and appending $r$ leaves to $y$ and $s$ leaves to $x$. It is worth noting that all trees of diameter 3 are double stars and vice versa.

\begin{figure}[ht]
\begin{tikzpicture}[scale=.75]
	\begin{pgfonlayer}{nodelayer}
		\node [style=blackoutline, label={above: $y$}] (0) at (-3.5, 0.5) {};
		\node [style=blackoutline, label={above: $x$}] (1) at (1.5, 0.5) {};
		\node [style=blackoutline, label={left: $y_2$}] (2) at (-7.75, 1.25) {};
		\node [style=blackoutline, label={left: $y_1$}] (3) at (-7.75, 2.25) {};
		\node [style=blackoutline, label={left: $y_r$}] (4) at (-7.75, -1.00) {};
		\node [style=blackoutline, label={right: $x_1$}] (5) at (5.75, 2.25) {};
		\node [style=blackoutline, label={right: $x_2$}] (6) at (5.75, 1.25) {};
		\node [style=blackoutline, label={right: $x_s$}] (7) at (5.75, -1.00) {};
		\path (4) -- node[auto=false]{\vdots} (2);
		\path (7) -- node[auto=false]{\vdots} (6);
	\end{pgfonlayer}
	\begin{pgfonlayer}{edgelayer}
		\draw [style=bold edge] (0) to (1);
		\draw [style=bold edge] (2) to (0);
		\draw [style=bold edge] (3) to (0);
		\draw [style=bold edge] (4) to (0);
		\draw [style=bold edge] (1) to (5);
		\draw [style=bold edge] (1) to (6);
		\draw [style=bold edge] (1) to (7);
	\end{pgfonlayer}
\end{tikzpicture}
\caption{$DS_{r,s}$}
\end{figure}

\section{Previous Results}

In \cite{kmsv}, the authors introduce the rainbow Tur\'an problem, and further conjecture that $\exr(n, C_{2k})= O(n^{1+ 1/k})$. In 2013, Das, Lee, and Sudakov show $\exr(n, C_{2k}) = O(n^{1+ \frac{(1+ \epsilon_k)lnk}{k}})$ with $\epsilon_k~\rightarrow~0$ as $k~\rightarrow~\infty$ in \cite{evencycles-das}. O. Janzer proves the conjecture in its entirety and further extends the result to cover all theta graphs in \cite{evencycles-janzer}. A 2016 paper by Johnston, Palmer and Sarkar \cite{pathandforest} disproves a conjecture on paths from \cite{kmsv} showing that $\exr(n,P_l)\le\lceil \frac{3l-2}{2}\rceil n$, and provides exact results for some particular $l$. Shortly after \cite{path-erg} improved this bound to $\exr(n, P_l) < \left(\frac{9(l-1)}{7} +2\right)n$. A lower bound for the rainbow Tur\'{a}n number of paths is given in \cite{puckanddan} as $\exr(n, P_k) \geq \frac{k}{2}n+ O(1)$. The same paper provides bounds on the rainbow Tur\'{a}n number of caterpillars and brooms, as well as an exact result for a specific family or brooms.

\section{Summary Of Results}

From \cite{kmsv} we have that the rainbow Tur\'{a}n number satisfies the inequality $\ex(n,F)~\leq~\exr(n,F)~\leq~\ex(n,H)~+~o(n^2)$. In conjunction with the Erd\H{o}s-Stone theorem, this result gives the order of magnitude for the rainbow Tur\'{a}n number of bipartite graphs. Therefore we restrict our consideration of rainbow Tur\'{a}n numbers to bipartite graphs. In particular, we present results for the family of double stars as defined in Section 1. Additionally, we present results on $k$-unique Tur\'{a}n numbers. These are of particular interest due to the fact that we have the following chain of inequalities for any graph $F$:
\begin{equation*}
    \ex(n,F)=\ex_0(n,F)\leq\ex_1(n,F)\leq\cdots\leq\ex_{\|H\|}(n,F)=\exr(n,F)
\end{equation*}
Part of the motivation for the $k$-unique Tur\'{a}n number is that they could allow us to improve existing lower bounds for the rainbow Tur\'{a}n numbers by finding $k$-unique Tur\'{a}n numbers and increasing the value of $k$. 

In order to consider $k$-unique Tur\'{a}n numbers, we begin by determining which exactly-$k$-unique values are possible for particular graphs. With that in mind, we define the {\bf $k$-spectrum} of a graph $F$, as the set of $k$ for which $F$ admits an exactly $k$-unique-coloring. We say that a graph has {\bf{full spectrum}} if $\spec(F)=\{0,1,\ldots,\|F\|-2,\|F\|\}$.\footnote{We note that having $(\|F\|-1)$ distinctly colored edges is impossible. This leaves one edge whose color can neither be unique nor match the color of any other edges.} 

Theorem \ref{thm:compspec} gives a criteria for some graphs which have a complete $k$-spectrum. Then, in Corollary \ref{cor:cyclespec} and Corollary \ref{cor:pathspec} we present a complete description of the $k$-spectrum of cycles and paths respectively. We conclude our focus on $k$-spectrum with Lemma \ref{lem:speclemma} which describes the $k$-spectrum for all double stars. Theorem \ref{thm:exkds} presents an upper bound for $\ex_k(n, DS_{r,s})$ when $DS_{r,s}$ admits an exactly-$k$-unique coloring. The following Lemma \ref{lem:exnkds} extends the result to cover $\ex_k(n, DS_{r,s})$ when $DS_{r,s}$ does not admit an exact $k$-unique coloring. 

The upper bound in Lemma \ref{thm:exkds} is achieved through a combination of a reduction method inspired by the method in \cite{kmsv} and the Erd\H{o}s-S\'{o}s Conjecture \cite{esconj}. Proposed in 1963 by Vera S\'{o}s, the Erd\H{o}s-S\'{o}s Conjecture is based on the observation that the traditional Tur\'{a}n number is the same for paths and stars. 
\begin{conj}
For any tree $T$ with $t$ edges, $\ex(n,T) \leq \frac{(t-1)n}{2}$.
\end{conj}
In the early 2000's, Ajtai, Koml\'{o}s, Simonovits, and Szemer\'{e}di announced a proof of the conjecture; this has not yet appeared. There have been specific cases of the Erd\H{o}s-S\'{o}s Conjecture for which results have been published. Of particular note is a result by McLennan \cite{mclennan2005erdHos}] which proves the Erd\H{o}s-S\'{o}s conjecture for trees with diameter at most 4. The double star (along with short brooms and caterpillars) fall under the above diameter constraint. 

Theorem \ref{thm:ds} uses the same methods to give an upper and lower bound bound on $\exr(n, DS_{r,s})$. In Theorems \ref{thm:exrbroom}, \ref{thm:dsk61}, and \ref{thm:exds2} we use the previous upper bound, but construct graphs which provide better lower bounds for the rainbow Tur\'{a}n numbers of $DS_{2,2}$ and $DS_{1,2k+1}$. Our result in Theorem \ref{thm:exrbroom} matches a result of Johnston and Rombach in \cite{puckanddan}, however we use a different method to achieve the matching upper bound. 

In Theorem \ref{thm:caterpillars}, \ref{thm:binary}, and \ref{thm:kary} we extend our reduction method to cover caterpillars and $k$-ary trees. In general, these trees do not fall under the result from \cite{mclennan2005erdHos}. In this case, we assume the Erd\H{o}s-S\'{o}s Conjecture is true for all trees and proceed by creating an algorithm to augment the caterpillars and $k$-ary trees respectively in order to reduce the rainbow Tur\'{a}n problem into a traditional Tur\'{a}n problem. 


\section{Proofs}

\begin{thm}
 Let $F$ be a graph with a proper edge coloring with color classes $L_1,\ldots,L_r$, ordered so that $|L_1|\ge|L_2|\ge\cdots\ge|L_r|$. If $|L_1|\ge 3$ and $|L_r|\ge 2$, then $F$ has full spectrum, that is $\spec(F)= \{0, 1, \ldots, \|F\|-2, \|F\|\}$ \label{thm:compspec}
\end{thm}

\begin{proof}
 Let $F$ be a graph with proper edge coloring with color classes $L_1,\ldots,L_r$ with $|L_i|=l_i$, ordered so that $l_1\ge l_2\ge\cdots\ge l_r$ and let $l_1\ge 3$ and $l_r\ge 2$. Label the edges of $F$ by $e_{i,j}$ with $1 \leq j \leq l_i$ such that $L_i=\{e_{i,1}, \ldots, e_{i,l_i}\}$. Define a new proper edge coloring $\phi_R:E(F) \rightarrow \|F\|$ where each edge is assigned a distinct color, with the restriction that $\phi_R(e_{i,1})=i$. 

Under $\phi_R$, each edge of $F$ has a distinct color, so $\|F\| \in \spec(F)$. Recoloring $e_{1,2}$ to color 1 gives a $(\|F\|-2)$-unique coloring since neither $e_{1,1}$ nor $e_{1,2}$ are assigned a distinct color. Recoloring the additional edges $e_{1, j}$ to color 1 will decrease the number of distinctly colored edges by one because $e_{1,1}$ is no longer the only edge in its color class. Doing so sequentially gives $\kappa$-unique colorings for $(\|F\|-3) \leq \kappa \leq (\|F\|-l_1)$, and so each of those values must belong in $\spec(F)$.

We know $l_1 \geq 3$, but it is possible that $l_i =2$ for all $i > 1$. Additionally, if we simply recolor $e_{i,2}$ with color $i$, the number of distinctly colored edges will decrease by two as it did for $e_{1,2}$ and $L_1$. The following color-switching operation addresses both of these problems simultaneously. 

Start with the coloring in which each of the $l_1$ edges $e_{1,j}$ are color 1 and all remaining edges in $F$ are colored by $\phi_R$. Then simultaneously recolor $e_{1,l_1}$ to color $\phi_R(e_{1,l_1})$ and $e_{2, 2}$ to color 2. By assigning $e_{1,l_1}$ a distinct color for this step, a $(\|F\|-l_1 -1)$-unique coloring is constructed. To construct a $(\|F\|-l_1-2)$-unique coloring, change the color of $e_{1,l_1}$ back to color 1. If $l_2 >2$, each remaining edge $e_{2,j}$ can be reassigned to color 2 one at a time.  As before, this process constructs $\kappa$-unique colorings of $F$ with $(\|F\|-l_1-3) \leq \kappa \leq (\|F\|-l_1-l_2)$. If $l_2 = 2$, there are no more edges $e_{2,j}$.

Repeating this color switching process using $e_{1, l_1}$ and $e_{i,2}$ for all $2 \leq i \leq r$ constructs a set of colorings that show $\spec(F)$ is complete.
\end{proof}

\begin{cor} Every cycle $C_k$ with $k\geq 6$ has full spectrum. For cycles with less than $6$ vertices, we have $\spec(C_5)=\{1,3,5\}$, $\spec(C_4)=\{0,2,4\}$, and $\spec(C_3)=\{3\}$. \label{cor:cyclespec}
\end{cor}
\begin{proof}
The result for $C_k$ with $k \geq 6$ follows directly from Theorem \ref{thm:compspec}.  For $C_k$ with $k < 6$, we proceed by cases.  Every proper edge coloring of $C_3$ is equivalent to the rainbow coloring, so $\spec(C_3)=\{3\}$. 

 Every proper edge coloring of $C_4$ which uses two colors is equivalent. Further, each color class of has cardinality two, so no edges are assigned a unique color and so $0 \in \spec(C_4)$. There are two non-equivalent proper edge colorings of $C_4$ which use exactly three colors. In both there is a single color class with cardinality two and two color classes with cardinality one. For every proper edge coloring on $C_4$ with three colors, there are exactly two edges assigned a distinct color, and $2 \in \spec(C_4)$. Finally, we consider the edge coloring in which each edge is assigned a distinct color. There are only 4 edges, so we have $4 \in \spec(C_4)$. Then $\spec(C_4) = \{0,2,4\}$.

$C_5$ has a chromatic index of $3$, so we need not consider edge colorings using fewer than $3$ colors. There are 5 non-equivalent proper edge colorings of $C_5$ using three colors, all with two color classes of cardinality two and one color class with cardinality one. These colorings each assign a single edge to a unique color, so $1 \in \spec(C_5)$. There are 5 non-equivalent proper edge colorings using precisely 4 colors. In each of these there is a single color class with cardinality two and three color classes with cardinality one. Thus $3 \in \spec(C_5)$. The coloring which assign each edge a distinct color uses 5 colors, and thus $5 \in \spec(C_5)$. then we find that $\spec(C_5) = \{1,3,5\}$.
\end{proof}

\begin{cor}
Every path $P_k$ of length $k$ with $k \geq 5$ has full spectrum. For paths of length less than 5, we have $\spec(P_4) = \{0,2,4\}$, $\spec(P_3)= \{1, 3\}$, $\spec(P_2)=\{2\}$, and $\spec(P_1)=\{1\}$. \label{cor:pathspec}
\end{cor} 
\begin{proof}
The result for $P_k$ with $k \geq 5$ follows directly from Theorem \ref{thm:compspec}. For $P_k$ with $k< 5$ we proceed by cases. $P_1$ and $P_2$ each have a single proper edge coloring (up to a renaming of the colors), and they both assign each edge a distinct color, so $\spec(P_1)=\{1\}$ and $\spec(P_2)= \{2\}$.

For $P_3$, there are exactly two non-equivalent proper edge colorings. The edge coloring of $P_3$ which two colors has one color class of cardinality two and one color class with cardinality one. Thus $1 \in \spec(P_3)$. Now we consider the edge coloring which assigns each edge to a distinct color. This gives $3 \in \spec(P_3)$, and thus $\spec(P_3)=\{1,3\}$.

 There are three kinds of non-equivalent proper edge colorings of $P_4$. Every proper edge coloring of $P_4$ which uses exactly two colors is equivalent, all with two color classes with cardinality two. In each, no edge is assigned to a distinct color so $0 \in \spec(P_4)$. For $P_4$ with precisely three colors, there are two non-equivalent proper edge colorings. In each there are two color classes with cardinality one and one color class with cardinality two. Thus $2 \in \spec(P_4)$. Finally, the edge coloring which assigns each edge a distinct color gives $4 \in \spec(P_4)$, and we have $\spec(P_4) = \{0,2,4\}$.
\end{proof}

\begin{lem} For a double star $\ds$ with $r \leq s$ and $j=s-r+1$, we have \begin{equation*}\spec(\ds) = \{j+2l: 0\leq l \leq r\}.\end{equation*}
\label{lem:speclemma}\end{lem}
\begin{proof} 
Since $yx$ is a dominating edge, its color cannot appear on any other edge of $\ds$. This means that in any proper edge coloring, any repeated colors must appear as pairs of pendants, one each from $y$ and $x$. Then no color can appear more than twice. In particular the number of repeated colors is at most $r$.  If every color which appears on $y$ also appears on a pendant from $x$, the remaining $s-r$ pendants from $x$ are distinctly colored (giving us $j=s-r+1$ distinctly colored edges from these pendants and the dominating edge). Now, by recoloring $l$ pendants of $y$ with colors not appearing elsewhere in the graph, we obtain a coloring with $j+2l$ distinctly colored edges since each $y$ pendant and its partner $x$ pendant now have their own distinct colors (since $y$ has only $r$ pendants, we require that $l\leq r$).
\end{proof}

\begin{thm}
For a double star $\ds$ with $r \leq s$ which admits an exactly $(j+2l)$-unique-coloring, we have \begin{equation*}\frac{(s+l-1)n}{2} +o(n) \leq \ex_{j+2l}(n, \ds) \leq \frac{(r+s+l)n}{2}\end{equation*} where $j=s-r+1$ and $0 \leq l \leq r$. \label{thm:exkds}
\end{thm} 
\begin{proof}
From Lemma \ref{lem:speclemma}, we have $DS_{r,s}$ colored with $s+l-1$ colors is $(j+2l)$-unique. By Vizing's theorem, we know that a graph $G$ with $\Delta(G)=s+l-1$ can be properly edge colored with $s+l$ colors. Let $G$ be a $(s+r-1)$-regular graph, or an almost $(s+r-1)$-regular graph if $(s+r-1)n$ is not even. Then $G$ can be properly colored so there is no rainbow $DS_{r,s}$ subgraph. This gives $\frac{(s+l-1)n}{2} +o(n) \leq \ex_{j+2l}(n, DS_{r,s})$. 

Consider any proper edge coloring of $DS_{r,s+l}$. In this augmented graph we call the endpoints of the dominating edge $y$ and $x'$, where $d(x') = s+l$. As in the proof of the previous lemma, any repeated colors in $DS_{r,s+l}$ must appear as pairs of pendants from $y$ and $x'$.  Thus there remain $s+l-r$ distinctly colored children of $x'$; along with $l$ many of $y$'s children and $yx'$, we have $(s+l-r)+l+1$ which is $j+2l$ distinctly colored edges. Applying the bound from the Erd\H{o}s-S\'{o}s conjecture gives $\ex(n, DS_{r,s+l}) \leq \frac{(r+s+l)n}{2}$. Thus we have $\ex_{j+2l}(n,\ds) \leq \ex(n, DS_{r, s+l}) \leq \frac{(r+s+l)n}{2}$.
\end{proof}

\begin{lem}
Let $k$ be given and, and set $k'$ to be the smallest $k'\ge k$ for which $k'\in\spec(F)$.  Then $\ex_{k'}(n,F)=\ex_{k}(n,F)$, since no exactly $\kappa$-unique $F$ exists for $k<\kappa<k'$. \label{lem:exnkds}
\end{lem}
\begin{proof} 
Consider the extremal graph $G$ which can be colored to avoid a $k$-unique copy of $F$.  Since no exactly $\kappa$-unique coloring is possible with $k~<~\kappa~<~k'$ (since such values of $\kappa$ are not in the spectrum), this $k$-unique copy must have at least $k'$ unique colors.
\end{proof}

\begin{thm}
For all double stars $\ds$ with $r \leq s$, we have 
\begin{equation*}
 \frac{(s+r-1)n}{2} +o(1)\leq \exr(n, \ds) \leq \frac{(s+2r)n}{2}.
\end{equation*}\label{thm:ds} \end{thm} 
\begin{proof} 
It follows from Theorem \ref{thm:exkds} that any properly edge colored $DS_{r, s+r}$ must contain a rainbow $DS_{r,s}$. Therefore we have the following upper bound: $\exr(n, DS_{r,s}) \leq \ex(n,DS_{r, s+r}) \leq \frac{(s+2r)n}{2}$.
The lower bound is shown by a simple application of Vizing's theorem, as in Theorem \ref{thm:exkds}. In order for every proper edge coloring of some graph $G$ to force a rainbow copy of $\ds$, it is necessary that every proper coloring of $G$ uses at least $s+r+1$ colors. By Vizing's theorem, a graph with $\Delta(G)+1=s+r$ can be properly colored to avoid a rainbow $\ds$ subgraph. This gives a lower bound of $\exr(n, \ds) \geq \frac{(s+r-1)n}{2} +o(1)$. 
\end{proof}

\begin{thm} Let $\phi:E(K_6) \rightarrow [5]$ be an edge coloring of $K_6$ where every vertex is adjacent to exactly one edge of each color. Then $\phi$ is rainbow $DS_{2,2}$-free. \label{thm:dsk61}
\end{thm}
\begin{proof} 
 Let $\phi:E(K_6) \rightarrow [5]$ be an edge coloring of $K_6$ where every vertex is adjacent to exactly one edge of each color, and suppose that $K_6$ with $\phi$ contains a rainbow $DS_{2,2}$. Without loss of generality, choose any vertex to act as $y$ in a $DS_{2,2}$. Under $\phi$, this vertex must be adjacent to five different colored edges. In the image below, the center vertex is $y$, and the solid black edges are assigned colors 1 through 5.

 Again, without loss of generality, we select a vertex adjacent to $y$ to be $x$ and choose which vertices will be their leaves as shown by the vertex labels in Figure \ref{fig:k6}. In order for this $DS_{2,2}$ to have each edge assigned to a distinct color, edges $xx_{1}$ and $xx_{2}$ (shown below with big dashes) must be assigned colors 5 and 2 respectively. The remaining edges adjacent to $x$, specifically $xy_1$ and $xy_2$ (shown below with small dashed edges) must be assigned colors 4 and 3 respectively. 
 \begin{figure}[ht]
\begin{tikzpicture}[scale=.4]
	\begin{pgfonlayer}{nodelayer}
		\node [style=blackoutline, label={above:$x$}] (0) at (0, 4.75) {};
		\node [style=blackoutline, label={below:$y_1$}] (1) at (-3.25, -4) {};
		\node [style=blackoutline, label={below:$y_2$}] (2) at (3.25, -4) {};
		\node [style=blackoutline, label={right:$x_2$}] (3) at (4.75, 1.5) {};
		\node [style=blackoutline, label={left:$x_1$}] (4) at (-4.75, 1.5) {};
		\node [style=blackoutline, label={below:$y$}] (5) at (0, 0) {};
		\node [style=none] (6) at (0.1, 0.75) {1};
		\node [style=none] (7) at (-1, 0.25) {2};
		\node [style=none] (8) at (-0.80, -1.00) {3};
		\node [style=none] (9) at (0.80, -1.00) {4};
		\node [style=none] (10) at (1, .25) {5};
	\end{pgfonlayer}
	\begin{pgfonlayer}{edgelayer}
		\draw [style=bold edge] (0) to (5);
		\draw [style=bold edge] (5) to (3);
		\draw [style=bold edge] (5) to (2);
		\draw [style=bold edge] (5) to (1);
		\draw [style=bold edge] (4) to (5);
		\draw [style=black dash 1] (0) to (4);
		\draw [style=black dash 1] (0) to (3);
		\draw [style=little dots] (0) to (1);
		\draw [style=little dots] (0) to (2);
		\draw [style=grey edge 1] (4) to (3);
		\draw [style=grey edge 1] (4) to (2);
		\draw [style=grey edge 1] (4) to (1);
	\end{pgfonlayer}
\end{tikzpicture}
    \caption{$DS_{2,2}$ embedded in $K_6$.}
    \label{fig:k6}
\end{figure}
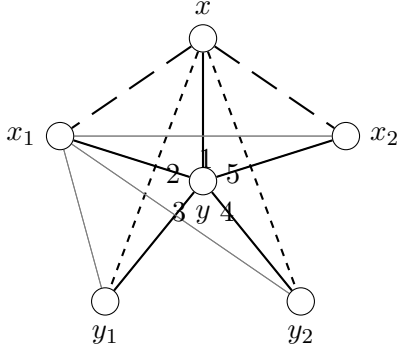    

Then we consider how colors can be assigned to the uncolored edges incident to vertex $x_1$ (shown in grey). The three colors 1, 3, and 4 have not met vertex $x_1$. However neither edge $x_1y_1$ nor $x_1y_2$ may be assigned color 3 or 4. This is a contradiction, since in $\phi$ every vertex is incident to an edge of every color. Thus there is no rainbow $DS_{2,2}$ in $K_6$ with edge coloring $\phi$. 
\end{proof}

\begin{thm}
  Any proper edge coloring of $K_6$  which is not rainbow contains an exactly-$3$-unique $DS_{2,2}$.  \label{thm:dsk62}
\end{thm}
\begin{proof}Consider any proper edge coloring $\phi~:~E(K_6)~ \rightarrow~\mathbb{N}$ which is not rainbow. Then two edges must assigned to the same color. Label the vertices such that  $v_1v_2$ and $v_4v_5$ are one of these pairs as in Figure \ref{fig:k6v2}.  

\begin{figure}[ht]
\begin{tikzpicture}[scale=.4]
	\begin{pgfonlayer}{nodelayer}
		\node [style=blackoutline, label={above:$v_6$}] (0) at (0, 5.25) {};
		\node [style=blackoutline, label={below:$v_3$}] (1) at (0, -5.25) {};
		\node [style=blackoutline, label={left:$v_1$}] (2) at (-5, 2.5) {};
		\node [style=blackoutline, label={left:$v_2$}] (3) at (-5, -2.5) {};
		\node [style=blackoutline, label={right:$v_5$}] (4) at (5, 2.5) {};
		\node [style=blackoutline, label={right:$v_4$}] (5) at (5, -2.5) {};
	\end{pgfonlayer}
	\begin{pgfonlayer}{edgelayer}
		\draw [style=bold edge] (2) to (5);
		\draw [style=black dash 1] (2) to (0);
		\draw [style=black dash 1] (0) to (5);
		\draw [style=black dash 1] (5) to (1);
		\draw [style=black dash 1] (1) to (2);
		\draw [style=grey edge 1] (5) to (4);
		\draw [style=grey edge 1] (3) to (2);
	\end{pgfonlayer}
\end{tikzpicture}
    \caption{$DS_{2,2}$ embedded in $K_6$ with repeated edge color.}
    \label{fig:k6v2}
\end{figure}
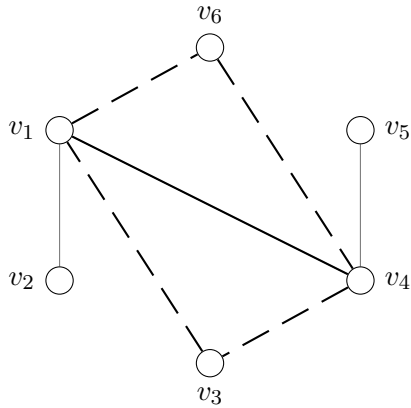

Let $v_1v_4$ be the dominating edge of a double star. From the two remaining vertices $v_3$ and $v_6$, one must be a leaf incident to $v_1$ and one must be a leaf incident to $v_4$. Note the $C_4$: $v_1v_6v_4v_3$. Regardless of how $\phi$ assigns colors to these four edges, we are able to choose our two leaves so the edges have distinct colors, thus creating an exactly-$3$-unique $DS_{2,2}$ subgraph.

\end{proof}

\begin{thm}
For the double star $DS_{2,2}$, we have the following inequality: \begin{equation*}
    \frac{5n}{2} \leq \exr(n, DS_{2,2}) \leq \frac{6n}{2}.
\end{equation*} \label{thm:exds2}
\end{thm}

\begin{proof}
The lower bound follows from Theorem \ref{thm:dsk61}. Since $K_6$ is 5-regular and can be properly colored with no rainbow $DS_{2,2}$, we know a graph must contain more than $\frac{5n}{2}$ edges in order to guarantee that every proper edge coloring contains a rainbow $DS_{2,2}$. For the upper bound, we apply the result from Theorem \ref{thm:ds}. Since $r=s=2$, we find $\exr(n,DS_{2,2}) \leq \frac{6n}{2}$.
\end{proof}

\begin{thm}
For $DS_{1,2s+1}$ for any $s\in\mathbb{N}$,
\begin{equation*}
\exr(n, DS_{1, 2s+1}) = \frac{(2s+3)n}{2} + o(1).\end{equation*} \label{thm:exrbroom}
\end{thm}
\begin{proof} The upper bound is an application of the result from Theorem \ref{thm:ds} which gives $\exr(n, DS_{1,2s+1}) \leq \frac{(2s+3)n}{2}$. We prove the lower bound by showing that it is possible to properly edge color $K_{2s+4}$ in a way that avoids a rainbow $DS_{1, 2s+1}$. Note that $|K_{2s+4}|=|DS_{1,2s+1}|$. Further, since $2s+4$ is even, $K_{2s+4}$ admits an edge coloring with $2s+3$ colors in which each vertex meets each color exactly once. We call this edge coloring $\phi$. 
Assume $K_{2s+4}$ with $\phi$ contains a rainbow $DS_{1, 2s+1}$. Then $DS_{1, 2s+1}$ must be embedded in the graph as in Figure \ref{fig:ds12k1}, with the edge colors as labeled (without loss of generality).

\begin{figure}[ht]
    \begin{tikzpicture}[scale=.4]
	\begin{pgfonlayer}{nodelayer}
		\node [style=blackoutline] (0) at (0, 5.25) {};
		\node [style=blackoutline] (1) at (0, -5.25) {};
		\node [style=blackoutline, label={left:$x$}] (2) at (-5, 2.5) {};
		\node [style=blackoutline, label={left:$y$}] (3) at (-5, -2.5) {};
		\node [style=blackoutline] (4) at (5, 2.5) {};
		\node [style=blackoutline] (5) at (5, -2.5) {};
		\node [style=none] (6) at (-5.25, 0) {1};
		\node [style=none] (7) at (-3, -3.75) {2};
		\node [style=none] (8) at (-2, 4.25) {3};
		\node [style=none] (9) at (1.00, 2.5) {4};
		\node [style=none] (10) at (1, -0.5) {$2k+3$};
		\node [style=none] (11) at (-0.5, 1.5) {};
		\node [style=none] (12) at (-0.5, 1.5) {$\vdots$};
	\end{pgfonlayer}
	\begin{pgfonlayer}{edgelayer}
		\draw [style=bold edge] (2) to (3);
		\draw [style=bold edge] (3) to (1);
		\draw [style=bold edge] (2) to (0);
		\draw [style=bold edge] (2) to (4);
		\draw [style=bold edge] (2) to (5);
		\draw [style=black dash 1] (2) to (1);
	\end{pgfonlayer}
\end{tikzpicture}
    \caption{$DS_{1, s2+1}$ embedded in $K_{2s+4}$.}
    \label{fig:ds12k1}
\end{figure}
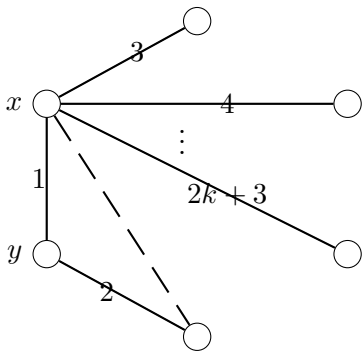

However, the dotted edge can not be assigned any of the $2s+3$ colors and a $(2s+4)^{th}$ color must be used. This is a contradiction to $\phi$. This gives $\exr(n, DS_{1, 2s+1}) \geq \frac{(2s+3)n}{2} + o(1)$ with equality when $(2s+3)$ divides $n$.
\end{proof}

A broom $B_{k,r}$ is a path on $k$ vertices with $r$ leaves appended to one endpoint. Every graph in Theorem \ref{thm:exrbroom} is also a broom. However, there is a more general family of graphs, caterpillars $C_{c_1,\ldots, c_k}$, consisting a path on $k$ vertices labeled $x_1, \ldots x_k$ with $c_i$ pendants attached at vertex $x_i$. This means every broom $B_{k,r}$ is also a caterpillar $C_{0, \ldots, 0, r}$, where $r$ appears in the $k^{th}$ position. With this in mind, the following theorem also applies to all brooms, although Theorem \ref{thm:exrbroom} gives a better result when $k = 3$. 

\begin{figure}[ht]
   \begin{tikzpicture}[scale=.42]
	\begin{pgfonlayer}{nodelayer}
		\node [style=blackoutline, label={above:$x_2$}] (0) at (-12, 0) {};
		\node [style=blackoutline, label={above:$x_1$}] (1) at (-15, 0) {};
		\node [style=blackoutline, label={above:$x_k$}] (2) at (-6, 0) {};
		\node [style=blackoutline, label={above:$x_3$}] (3) at (-9, 0) {};
		\node [style=blackoutline, label={right:$y_2$}] (4) at (-3, 0) {};
		\node [style=blackoutline, label={right:$y_1$}] (5) at (-3, 2.5) {};
		\node [style=blackoutline, label={right:$y_r$}] (6) at (-3, -2.75) {};
		\node [style=none] (7) at (-3, -1.5) {};
		\node [style=none] (8) at (-3, -1.5) {$\vdots$};
		\node [style=none] (9) at (-7.5, 0) {$\ldots$};
		\node [style=blackoutline, label={above:$x_2$}] (10) at (6, 0) {};
		\node [style=blackoutline, label={above:$x_1$}] (11) at (2, 0) {};
		\node [style=blackoutline, label={above:$x_k$}] (12) at (14, 0) {};
		\node [style=blackoutline, label={above:$x_3$}] (13) at (10, 0) {};
		\node [style=blackoutline, label={below:$y_{k,1}$}] (15) at (13, -2) {};
		\node [style=blackoutline, label={below:$y_{k,c_k}$}] (16) at (15, -2) {};
		\node [style=none] (19) at (12, 0) {$\ldots$};
		\node [style=none] (20) at (14, -2) {$\ldots$};
		\node [style=blackoutline, label={below:$y_{1,1}$}] (21) at (1, -2) {};
		\node [style=blackoutline, label={below:$y_{1,c_1}$}] (23) at (3, -2) {};
		\node [style=blackoutline, label={below:$y_{2,1}$}] (24) at (5, -2) {};
		\node [style=blackoutline, label={below:$y_{2,c_2}$}] (25) at (7, -2) {};
		\node [style=blackoutline, label={below:$y_{3,1}$}] (26) at (9, -2) {};
		\node [style=blackoutline, label={below:$y_{3,c_3}$}] (27) at (11, -2) {};
		\node [style=none] (28) at (10, -2) {$\ldots$};
		\node [style=none] (29) at (6, -2) {$\ldots$};
		\node [style=none] (30) at (2, -2) {$\ldots$};
	\end{pgfonlayer}
	\begin{pgfonlayer}{edgelayer}
		\draw [style=bold edge] (1) to (0);
		\draw [style=bold edge] (0) to (3);
		\draw [style=bold edge] (2) to (5);
		\draw [style=bold edge] (2) to (4);
		\draw [style=bold edge, in=135, out=-45] (2) to (6);
		\draw [style=bold edge] (11) to (10);
		\draw [style=bold edge] (10) to (13);
		\draw [style=bold edge] (12) to (15);
		\draw [style=bold edge] (12) to (16);
		\draw [style=bold edge] (13) to (26);
		\draw [style=bold edge] (13) to (27);
		\draw [style=bold edge] (10) to (25);
		\draw [style=bold edge] (10) to (24);
		\draw [style=bold edge] (11) to (21);
		\draw [style=bold edge] (11) to (23);
	\end{pgfonlayer}
\end{tikzpicture}

    \caption{Broom $B_{k,r}$ (left) and caterpillar $C_{c_1, \ldots c_k}$ (right)}
    \label{fig:brromandcat}
\end{figure}
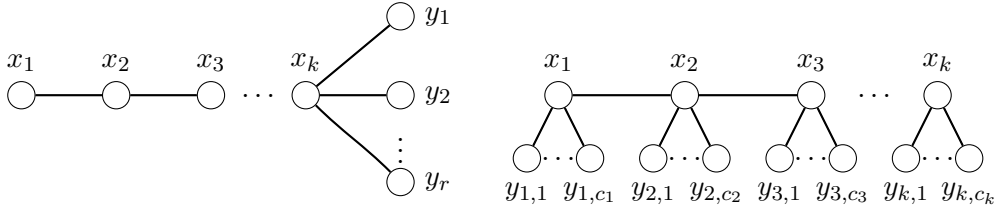

\begin{thm} Let $C_{c_1,\ldots, c_k}$ be a path on $k$ vertices labeled $x_1, \ldots x_k$ with $c_i$ pendants attached at vertex $x_i$. Then \[\exr(n, C_{c_1, \ldots c_k}) \leq [3c_1 +2c_2 + c_3+ 3 + 
\sum^k_{j=3}P_j(L_j +1)]\frac{n}{2}\] where $P_j = P_{j-1} \sum_{i=4}^{j-2}(c_i + 1)$ and $L_j= j-1 + \sum_{i=1}^{j}c_i$. \label{thm:caterpillars}
\end{thm}
\begin{proof}
We label the pendant vertices of $C_{c_1,\ldots, c_k}$ such that the children of $x_i$ are $v_{i,j}$ for $1 \leq j \leq c_i$. We will construct an augmented caterpillar $C_{c_1',\ldots,c'_k}$ which contains a rainbow $C_{c_1,\ldots, c_k}$ under every proper edge coloring. For $C_{c_1,\ldots, c_k}$ with $k = 3$, the two path edges will be assigned distinct colors under any proper edge coloring. The pendants adjacent to $x_1$ will all be assigned distinct colors as well, and by adding an extra pendant, we can ensure that there will be at least $c_1$ edges which do not share a color with either of the backbone edge. Then $c'_1 = c_1 +1$ in $C_{c'_1,\ldots, c'_k}$. The pendant edges adjacent to $x_2$ will all be assigned colors distinct from each other and distinct from both backbone edges. In order to ensure that there are $c_2$ of them which do not share a color with the pendants adjacent to $x_1$ we add another $c_1$ pendants. Thus $c'_2 = c_2 + c_1$. Finally we consider the pendant edges adjacent to $x_3$. There are $c_3$ of them in $C_{c_1,\ldots, c_k}$, and they will all be assigned colors distinct from each other and distinct from the edge $x_2x_3$. To guarantee that there are $c_3$ of them which do not share a color with the other backbone edge or any of the other pendant edges, we can simply add $c_1+c_2+1$ extra pendants. Then $c'_3 = c_1 + c_2 + c_3 + 1$. This new graph, $C_{c'_1,\ldots, c'_k}$, has $3c_1 + 2c_2 + c_3 + 5$ total edges. 

When $k= 3$, the augmented graph that must contain a rainbow $C_{c_1,c_2,c_3}$ subgraph is still a caterpillar. For higher values of $k$, this is no longer the case. As we generalize this method to cover larger $k$, we refer to the augmented graph of $C_{c_1, \ldots, c_k}$ as $C'_{c_1, \ldots, c_k}$. 

Consider $C_{c_1,\ldots, c_k}$ with $k=4$ as an augmentation of $C_{c_1,c_2,c_3}$. Then $C'_{c_1,\ldots, c_4}$ is also an augmentation of  $C_{c'_1,c'_2, c'_3}$. In order to ensure that there is an edge which can act as the image of $x_3x_4$ in the rainbow subgraph, add $c_1 + c_2 + 1$ more edges adjacent to $x_3$, for a total of $c_1+c_2+2$ non-pendant edges separate from $x_2x_3$. Each of these edges has $c_1+c_2+c_3+3$ pendant edges adjacent. This is large enough to ensure that as we find our rainbow subgraph we can choose an edge to serve as $x_3x_4$, along with the necessary $c_4$ pendants avoiding colors already in our subgraph (since there are at most $c_1+c_2+c_3+3$ such colors to avoid). 

Note that if we continue in this way for higher values of $k$, the next backbone branches will be adjacent to $x_{k-1}$ and not any of the $v_{k-1,j}$. With this in mind, the previous argument for $C_{c_1,\ldots, c_4}$ generalizes into the following parts: $P_j$ the number of backbone vertices which can act as $x_j$ in a rainbow subgraph, and $L_j$ the number of leaves adjacent to each of the $P_j$ parent vertices. The sum of all $c_i$ up to and including $j$ plus $(j-1)$ non-adjacent backbone edges is the number of leaves at level $j$. Then there must be  $L_j =j-1+ \sum_{i=1}^{j}c_i$ pendants adjacent to each parent vertex. Since the backbone branches at every new vertex, we multiply the number of previous parents by the sum of all previous children to determine the number of parent vertices at each level. Then $P_j = P_{j-1} \sum_{i=4}^{j-2}(c_i + 1)$. In this case we start branching at $i=4$ since we use the argument for $C_{c_1,c_2,c_3}$ before then. Thus the total number of edges in $C'_{c_1, \ldots, c_k}$ is $3c_1 + 2c_2 + c_3 + 5 + \left[\sum^k_{j=3}P_j(L_j +1)\right]$. Applying the bound from the Erd\H{o}s-S\'{o}s Conjecture gives the result stated in the theorem. 
\end{proof}

Let a $k$-ary tree be a rooted tree with every non-leaf vertex having exactly $k$ children, and let a perfect $k$-ary tree be a $k$-ary tree with all leaves at the same depth, or distance from the rooted vertex. Then a perfect binary tree with a depth of 2 has one vertex of degree 2, two vertices with degree 3, and 4 leaves. We use $T(k,d)$ to denote a perfect $k$-ary tree with depth $d$. Applying our reduction method to $T(k,d)$ gives the following results.

\begin{thm} 
For a perfect binary tree with depth $d$,
\begin{equation*}\exr(n, T(2,d)) \leq \Big(\sum^{d}_{j=2} [2 \prod^j_{i=2}(2^i -3)]+1\Big) \frac{n}{2}\end{equation*} \label{thm:binary}
\end{thm}
\begin{proof} Summing the $2^i$ vertices at each level of the tree, we find a total of $\sum_{i=0}^{d}2^i~=2^{d+1}-1$ vertices in $T(2,d)$ -- and thus $2^{d+1}-2$ edges. We label the vertices and edges in $T(2,d)$ according to the following system, an example of which is shown in Figure \ref{fig:binarytree}. The vertices are labeled $v_{i,j}$, for $1 \leq j \leq 2^i$, where $i$ is the distance to the root vertex, and with vertices $v_{i+1, 2j}$ and $v_{i+1, 2j-1}$ as the children of $v_{i,j}$. We will identify edges using the indices of the endpoint furthest from $v_{0,1}$; e.g. edge $v_{1,2}v_{2,4}$ would be labeled $e_{2,4}$. In the vein of the reduction method used earlier, we will construct a tree $T'(2,d)$ which contains a rainbow $T(2,d)$ subgraph under every proper edge coloring. For clarity, we do this first for $T(2,2)$ before generalizing our construction to all $T(2,d)$. 

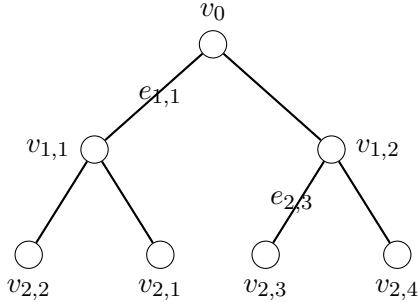
\begin{figure}[ht]
 \begin{tikzpicture}[scale=.7]
	\begin{pgfonlayer}{nodelayer}
		\node [style=blackoutline, label={above:$v_{0}$}] (0) at (0, 2) {};
		\node [style=blackoutline, label={left:$v_{1,1}$}] (1) at (-2.25, 0) {};
		\node [style=blackoutline, label={right:$v_{1,2}$}] (2) at (2.25, 0) {};
		\node [style=blackoutline, label={below:$v_{2,1}$}] (3) at (-1, -2) {};
		\node [style=blackoutline, label={below:$v_{2,2}$}] (4) at (-3.5, -2) {};
		\node [style=blackoutline, label={below:$v_{2,3}$}] (5) at (1, -2) {};
		\node [style=blackoutline, label={below:$v_{2,4}$}] (6) at (3.5, -2) {};
		\node [style=none] (7) at (-1, 1) {$e_{1,1}$};
		\node [style=none] (11) at (1.5, -1) {$e_{2,3}$};
	\end{pgfonlayer}
	\begin{pgfonlayer}{edgelayer}
		\draw [style=bold edge] (1) to (0);
		\draw [style=bold edge] (0) to (2);
		\draw [style=bold edge] (2) to (5);
		\draw [style=bold edge] (2) to (6);
		\draw [style=bold edge] (1) to (4);
		\draw [style=bold edge] (1) to (3);
	\end{pgfonlayer}
\end{tikzpicture}
    \caption{Perfect binary tree $T(2,2)$}
    \label{fig:binarytree}
\end{figure}

 In $T(2,2)$, $e_{1,1}$ and $e_{1,2}$ will be assigned distinct colors under any proper edge coloring. We need ensure that we can always find a pair of pendants for each of $v_{1,1}$ and $v_{1,2}$ whose colors are different from those already in our rainbow subgraph. Without loss of generality, we focus on the children of $v_{1,1}$. The two associated edges ($e_{2,1}$ and $e_{2,2}$) will always be assigned colors distinct from each other, and unique from $e_{1,1}$. In order to ensure there will be two edges with colors distinct from the rest of the leaves, we append two more leaves to the image vertex in $T'(2,2)$. The addition of one more leaf ensures that there is also an edge with a color unique from $e_{1,2}$. Then the total number of leaves adjacent to $v_{1,1}$ in $T'(2,2)$ is the total number of leaves in $T(2,2)$ plus the number of edges in $T(k, d-1)$ minus one. This gives $2^{d+1}-3$ leaves attached to each parent vertex.  Since there are two vertices at depth one, the total number of leaves in $T'(2,2)$ is 10. Including the edges $e_{1,1}$ and $e_{1,2}$, there are 12 edges in $T'(2,2)$. Applying the bound from the Erd\H{o}s-S\'{o}s Conjecture to $T'(2,2)$ gives $\exr(n, T(2,2)) \leq \ex(n,T'(2,2)) \leq \frac{10n}{2} $. Note this follows the form $(2(2^{d+1}-3))\frac{n}{2}$ .

The argument above generalizes neatly to cover all $T(2,d)$. 
Having $2^{d+1}~-~3$ leaf edges appended to each $v_{j-1,j}$ guarantees that even if  the colors of the other $2^d-2$ leaves and the colors of the $2^d-1$ non-adjacent edges in the the $T(2,d-1)$ subgraph are all repeated, there are still  two edges with distinct colors adjacent to $v_{d-1,j}$.  The number of \textbf{parents of leaves} in $T'(2,d)$ is the same as the \textbf{number of leaves} in $T'(2,d-1)$. Continuing in this way, we find the number of vertices at depth $d-2$ in $T'(2,d)$ is the same as the number of leaves in $T'(2,d-2)$. As in the previous example, at depth 1 there are two single edges which will always be assigned distinct colors. This gives $2 \prod^k_{i=2}(2^i -3)$ as the total number of leaves in $T'(k,d)$. The sum $(2+\sum^{d}_{j=2}(\text{number of leaves in } T'(2, j))$ provides the total number of edges in $T'(2,d)$. Applying the bound from the Erd\H{o}s-S\'{o}s Conjecture to this number gives $\exr(n, T(2,d)) \leq \ex(n, T'(2,d)) \leq (\sum^{d}_{j=2} [2 \prod^j_{i=2}(2^i -3)]+1) \frac{n}{2}$.
\end{proof}

\begin{thm}
For a perfect $k$-ary tree with depth $d$,
\begin{equation*}\exr(n, T(k,d)) \leq \Big(k - 1 + \sum^{d}_{j=2} \big[k \prod^j_{i=2}(k^i + \frac{k^i - 1}{k-1} - 2)\big]\Big) \frac{n}{2}\end{equation*} \label{thm:kary}
\end{thm}
\begin{proof}
Generalizing the geometric series in Theorem \ref{thm:binary} shows the number of vertices in $T(k,d)$ is $\frac{k^{d+1}-1}{k-1}$, and number of edges is $\frac{k^{d+1}}{k-1}-1$. Following the same argument as in Theorem \ref{thm:binary}, each vertex $v_{d-1,j}$ in $T'(k,d)$ must be adjacent to $k^d + \frac{k^d-1}{k-1}-2$ leaves to guarantee there are enough to choose $k$ many pendants with distinct colors in a rainbow $T(k,d)$ subgraph. The number of parents of leaves in $T'(k,d)$ is the same as the number of leaves in $T'(k,d-1)$. This process gives $k [\prod^j_{i=2}(k^i + \frac{k^i - 1}{k-1} - 2)]$ as the total number of leaves in $T'(k,d)$. Each vertex $v_{i,j}$ except for $v_{0,1}$ is associated with exactly one edge $e_{i,j}$. Then counting the number of vertices at each depth $i \geq 2$ gives $\sum^{d}_{j=2} [k \prod^j_{i=2}(k^i + \frac{k^i - 1}{k-1} - 2)]$, which is the same as the number of edges in $T'(k,d)-v_{0,1}$.  By adding $k$ more edges to account for those adjacent to vertex $v_{0, 1}$, we find the total number of edges in $T'(k,d)$. Applying the bound from the Erd\H{o}s-S\'{o}s Conjecture to this edge count gives an upper bound of $\exr(n, T(k,d)) \leq \ex(n, T'(k,d) \leq (k - 1 + \sum^{d}_{j=2} [k \prod^j_{i=2}(k^i + \frac{k^i - 1}{k-1} - 2)]) \frac{n}{2}$.
\end{proof}

\section{Conclusion}
As we see in Theorem \ref{thm:exds2} and Theorem \ref{thm:exrbroom}, using the reduction method in combination with the Erd\H{o}s-S\'{o}s Conjecture gives a good upper bound in cases where only one or few edges can be added to guarantee a rainbow subgraph. At first glance it seems as though the small diameter of such graphs results in good upper bounds, but consider $DS_{20, 20}$. In order to guarantee a rainbow subgraph, 20 edges must be appended to the original, and the upper bound achieved through the reduction method and Erd\H{o}s-S\'{o}s Conjecture is seemingly further than what we expect the true value to be. On the contrary for the graph $DS_{1,49}$, which has the same diameter and number of edges as $DS_{20,20}$, the methods here give matching upper and lower bounds on $\exr(n,DS_{1, 39})$. This implies a more complicated set of criteria for which this combination method provides good bounds; more investigation is needed. 

There are many trees not covered by the methods in this paper. Subdivided stars, and others with a simple branching structure should also lend themselves nicely to the reduction method as outlined in this paper. In particular, the Erd\H{o}s-S\'{o}s Conjecture has been proved for subdivided stars in \cite{esspiders}, and so these are a natural next step. 

In general, there may be ways to improve the reduction method to give better results for those trees with large size or diameter. For example, if we know that a host graph has a certain connectivity then how does that impact the rainbow subgraphs it may contain? Aside from just the connectivity, knowing more about the the properties/invariants of the extremal graphs for specified forbidden subgraphs may provide a method to reduce the number of added edges needed when implementing the reduction method. Further, know from Erd\H{o}s-Stone that augmenting trees such that the new graph is not a tree will cause a jump in the extremal number of $F'$, however such a method may be a useful implement if the forbidden subgraph $F$ already contains a cycle.

This paper only give a brief introduction into the $k$-unique Tur\'{a}n numbers. There is a lot more work that can be done here. Additionally, we define generalized $k$-unique Tur\'{a}n problems, as a natural extension of generalized rainbow Tur\'{a}n problems as defined in \cite{genrainbowturan}. What is the maximum number of copies of $F$ in a properly edge-colored graph on $n$ vertices without a $k$-unique copy of $F$.  It is important to note that using a reduction method to find upper bounds for the $k$-unique Tur\'{a}n numbers may provide useful information, but it will not allow us to find better lower bounds for the rainbow Tur\'{a}n number. Thus, if the goal is to improve the lower bounds for rainbow Tur\'{a}n numbers, the focus must be on finding lower bounds for the $k$-unique Tur\'{a}n numbers.

\end{document}